\documentclass[12pt]{amsart}
\usepackage{amsmath}
\usepackage{amsfonts}
\usepackage{fullpage}
\newfont{\Bbe}{msbm9}
\newfont{\Bb}{msbm10}

\begin{document}
\newcommand{\re}{\mbox{\Bb R}}
\newcommand{\na}{\mbox{\Bb N}}   
\newcommand{\ze}{\mbox{\Bb Z}}
\newcommand{\smallna}{\mbox{\Bbe N}}

\title[Connectedness and the action of group of isometries]{The role of
       connectedness in the structure and the action
       of group of isometries of locally compact metric spaces}
\author[Manoussos]{Antonios Manoussos}
\address{Antonios Manoussos\\
        123, Sapfous St.\\
        176 75 Kallithea\\
        Athens\\
    GREECE}
\email{amanou@cc.uoa.gr}
\author[Strantzalos]{Polychronis Strantzalos}
\address{Polychronis Strantzalos\\
    Dept.\ of Mathematics\\
    University of Athens\\
    Panepistimioupolis\\
    GR-157 84, Athens\\
    GREECE}
\email{pstrantz@math.uoa.gr}

\begin{abstract}
By proving that, if the quotient space $\Sigma (X)$ of the connected
components of the locally compact metric space $(X,d)$ is compact, then the
full group $I(X,d)$ of isometries of $X$ is closed in $C(X,X)$ with respect
to the pointwise convergence topology, i.e., that $I(X,d)$ coincides in this
case with its Ellis' semigroup, we complete the proof of the following:

\bigskip

\noindent {\bf Theorem}

\begin{description}
\item[(a)]  {\em If $\Sigma (X)$ is not compact, $I(X,d)$ need not be
locally compact, nor act properly on $X$. }

\item[(b)]  {\em If $\Sigma (X)$ is compact, $I(X,d)$ is locally compact but
need not act properly on $X$. }

\item[(c)]  {\em If, especially, $X$ is connected, the action $(I(X,d),X)$
is proper. }
\end{description}

\end{abstract}

\maketitle

\section{Introduction}

Transformation groups have been proven an effective tool to investigate the
structure of locally compact spaces. Since we know more about locally
compact groups than about locally compact spaces, it is reasonable to use
special properties of groups and their actions in order to understand the
structure of spaces on which they act. In this direction, we treat in this
note the following twofold question:

Let $(X,d)$ be a locally compact metric space and $I(X,d)$ its group of
isometries. When (a) is $I(X,d)$ locally compact (always with respect to the
pointwise convergence topology), and (b) does it act properly on $X$?

In the case where $X$ is connected the first question is answered
affirmatively in \cite{d-w} (cf. also \cite[Ch. I, Th. 4.7]{kob-nom}). The
non--connected case is investigated in \cite{str}. It turned out in this
case, that the compactness of $\Sigma (X)$, the quotient--space of the
connected components of $X$, is a topological property, independed of the
each time considered admissible metric, which ensures the local compactness
of $\overline{I(X,d)}$, the closure of $I(X,d)$ in $C(X,X)$, the space of
the continuous maps $X\rightarrow X$ endowed with the pointwise convergence
topology, for all admissible metrics $d$. The question whether $I(X,d)$ is
closed in $C(X,X)$ remained open. In this note we fill this gap (cf. also 
\cite{aman}), i.e., we show that in the case where $\Sigma (X)$ is compact $%
I(X,d)$ coincides with its Ellis' semigroup, completing the proof of the
following:

\bigskip

\noindent {\bf Theorem} {\em Let $(X,d)$ be a locally compact metric space.
Denote by $\Sigma (X)$ the space of the connected components of $X$, and by $%
I(X,d)$ its group of isometries. Then }

\begin{enumerate}
\item  {\em If $\Sigma (X)$ is not compact, $I(X,d)$ need not be locally
compact, nor act properly on $X$.}

\item  {\em If $\Sigma (X)$ is compact, then }

\begin{description}
\item[(a)]  {\em $I(X,d)$ is locally compact, }

\item[(b)]  {\em the action $(I(X,d),X)$ is not always proper, and }

\item[(c)]  {\em in the special case where $X$ is connected, the action $%
(I(X,d),X)$ is proper. }
\end{description}
\end{enumerate}

For the sake of completeness, we give short and slightly improved proofs of
already published by the authors partial results, crucial for a unified
proof of the above theorem. Our treatment is based on the sets $%
(x,V_{x})=\{g\in I(X,d):\,gx\in V_{x}\}$, where $V_{x}$ is a neighborhood of 
$x\in X$. These sets form a neighborhood subbases for the identity, for the
pointwise convergence topology, the natural topology of $I(X,d)$.

\section{Generalities}

\noindent {\large {\bf 2.1}} The following simple examples answer 1 and 2(b)
of the above theorem.

\bigskip

\noindent {\bf Example } Let $X={\Bbb {Z}}$ with the discrete metric.
Obviously $\Sigma (X)$ is not compact. It can be easily seen that $I(X,d)$
is the group of the bijections of ${\Bbb {Z}}$, which is not locally compact
with respect to the pointwise convergence topology, therefore it cannot act
properly on a locally compact space.

\bigskip

\noindent {\bf Example } Let $X=Y\cup \{(1,0)\}\subset {\Bbb {R}^{2}}$ where 
$Y=\{(0,y):\;y\in {\Bbb {R} \}}$, and $d=\min \{1,\delta \}$, where $\delta $
denotes the Euclidean metric. As we shall see in \S 3, by Theorem 3.7, $%
I(X,d)$ is locally compact; however the action $(I(X,d),X)$ is not proper,
because the isotropy group of $(1,0)$ is not compact, since it contains the
translations of $Y$. So, the action of $I(X,d)$ on $X$ is not proper, even
if $X$ has two components.

\bigskip

Since the sets $(x,V_{x})$ as above form a neighborhood subbases for the
identity for the pointwise convergence topology in $I(X,d)$, the following
condition is necessary for the local compactness of $I(X,d)$:

\begin{description}
\item[(a)]  There exist $x_{i}\in X$, $i=1,\ldots ,m$ such that $%
\bigcap_{i=1}^{m}(x_{i},V_{x_{i}})$ is relatively compact in $C(X,X)$.

This condition becomes also sufficient, if additionally the following
condition is satisfied.

\item[(b)]  $I(X,d)$ is closed in $C(X,X)$.
\end{description}

So, to prove the locall compactness of $I(X,d)$, we have to ensure that both
of the above conditions are satisfied.

\section{The local compactness of $I(X,d)$}

The following is crucial for the investigation of the conditions 2.1(a) and
(b):

\bigskip

\noindent {\large {\bf 3.1}} {\bf Lemma } {\em Let $(X,d)$ be a locally
compact metric space, $F\subseteq I(X,d)$, and} 
\[
K(F)=\{x\in X:\,F(x)=\{fx:f\in F\} \,\mbox {{\em is relatively compact}}\}. 
\]
\noindent {\em Then $K(F)$ is an open and closed subset of $X$.} \newline

\noindent {\bf Proof}. Let $x\in K(F)$ and $A$ be a relatively compact
neighborhood of $\overline{F(x)}$. Let $S(B,\varepsilon)$ denote the $%
\varepsilon$--neighborhood of $B\subseteq X$, and $\varepsilon$ be such that 
$S(F(x),\varepsilon)\subset A$. It is easily seen that $S(x,\varepsilon)%
\subset K(F)$, hence $K(F)$ is open.

On the other hand, if $K(F)\ni x_{n}\rightarrow x$ and $x_{n}\in S(x,\eta)$,
there are $f_{i}\in F$ such that 
\[
\overline{F(x_{n})}\subset \bigcup _{i=1}^{k}f_{i}(S(x,2\eta )). 
\]
\noindent Since 
\begin{eqnarray*}
d(fx,f_{i}x_{n}) &\leq &d(fx,fx_{n})+d(fx_{n},f_{i}x)+d(f_{i}x,f_{i}x_{n}) \\
&=&d(x,x_{n})+d(fx_{n},f_{i}x)+d(x,x_{n}),
\end{eqnarray*}
\noindent if $i\in \{1,\ldots ,k\}$ is such that $fx_{n}\in f_{i}(S(x,2\eta
)) $, then $d(fx,f_{i}x_{n})<4\eta .$ Hence 
\[
fx\in S(f_{i}x_{n},4\eta )=f_{i}(S(x_{n},4\eta ))\subseteq f_{i}(S(x,5\eta
)). 
\]
\noindent Therefore $F(x)\subset \bigcup _{i=1}^{k}f_{i}(S(x,5\eta ),$ which
means that $F(x)$ is relatively compact, if $S(x,5\eta )$ is.

\bigskip

\noindent {\large {\bf 3.2}} {\bf Remark} In the sequel we assume that $%
\Sigma (X)$ {\em is compact} in the quotient topology via the natural map $%
q:X\rightarrow \Sigma (X).$ Note that $\Sigma (X)$ is a $T_{1}$--space, and
need not be Hausdorff. Nevertheless

\bigskip \noindent {\em $X$ is separable, hence second countable; so
sequences are adequate in C(X,X).}

The proof is similar to the lengthy one in \cite{sie} (see also 
\cite[Appendix 2]{kob-nom}).

\bigskip

\noindent {\large {\bf 3.3}} {\bf Lemma } {\em Let $(X,d)$ be a locally
compact metric space with compact space of connected components $\Sigma (X)$%
. Then condition 2.1(a) is satisfied.} \newline

\noindent {\bf Proof.} Let $V_{x}$ be a relatively compact neighborhood of $%
x\in X$. Then 
\[
(x,V_{x})=\{g\in I(X,d):\,gx\in V_{x}\} 
\]
\noindent is a neighborhood of the identity in $I(X,d)$. By Lemma 3.1, $%
K((x,V_{x}))$ is not empty, open and contains whole components of $X$,
therefore $q(K((x,V_{x}))$) is an open subset of $\Sigma (X)$. Since $\Sigma
(X)$ is compact, there are $x_{i}$, $i=1,\ldots ,m$, such that the
corresponding $q(K((x_{i},V_{x_{i}}))$)' s cover $\Sigma (X)$. This means
that $X= \bigcup_{i=1}^{m}K((x_{i},V_{x_{i}}))$, i.e., the neighborhood $%
F=\bigcap_{i=1}^{m}(x_{i},V_{x_{i}})$ of the identity has the property: $F(x)
$ is relatively compact in $X$, for every $x\in X$, therefore, by Ascoli' s
theorem, $F$ is relatively compact in $C(X,X)$.

\bigskip

\noindent {\large {\bf 3.4}} Now we proceed to prove that under assumption
that $\Sigma (X)$ is compact, $I(X,d)$ is a closed subspace of $C(X,X)$.
Because of Remark 3.2, the elements $f$ of the boundary of $I(X,d)$ in $%
C(X,X)$ are limits of sequences $\{f_{n}\in I(X,d),n\in {\Bbb {N}\}}$.
Obviously, such an $f$ preserves $d$; so the question is when $f$ is
surjective. If $\Sigma (X)$ is not compact this is not always true:

\bigskip

\noindent {\bf Example} Let $X={\Bbb {Z}}$ with the discrete metric. If $%
f_{n}(z)=z$ for $-n<z<0$, $f(-n)=0$, and $f_{n}(z)=z+1$ otherwise then $%
f_{n}\rightarrow f$, where $f(z)=z$ for $z<0$, and $f(z)=z+1$ for $z\geq 0.$
Hence each $f_{n}$ is an isometry, but $f$ is not surjective since $0\notin
f({\Bbb {Z}).}$

\bigskip

\noindent {\large {\bf 3.5}} {\bf Lemma} {\em If $\Sigma (X)$ is compact and 
$I(X,d)\ni f_{n}\rightarrow f$, then $f(X)$ is open and closed in $X$.} 
\newline

\noindent {\bf Proof.} By Lemma 3.1, it suffices to show that $f(X)=K(F)$,
where $F=\{f_{n}^{-1},\,n\in {\Bbb {N}\}}$. Indeed, since $%
d(f_{n}(x),f(x))=d(x,f_{n}^{-1}(f(x))$, we have $f_{n}^{-1}(f(x))\rightarrow
x$, so (since $X$ is locally compact) $f(x)\in K(F)$, for every $x\in X.$
Now, if $y\in K(F)$, we may assume $f_{n_{k}}^{-1}(y)\rightarrow x$ for some 
$x\in X$, because $F(y)$ is relatively compact in $X$, hence $f(x)=y$.

\bigskip

\noindent {\large {\bf 3.6}} {\bf Propositon} {\em If $(X,d)$ is a locally
compact metric space, and $\Sigma (X)$ is compact, then $I(X,d)$ is closed
in $C(X,X)$.} \newline

\noindent {\bf Proof.} Let $I(X,d)\ni f_{n}\rightarrow f\in C(X,X)$. We
prove that $f$ is surjective. Let $y\in X$. We denote by $S_{x}$ the
connected component containing $x\in X$, and by $S_{n}$ the component of $%
f_{n}^{-1}y$. If $\{S_{n},n\in {\Bbb {N}\}}$ has any constant subnet $%
\{S_{n_{i}},\,i\in I\}$, then $S_{n_{i}}=S_{0}$, for some $S_{0} \in \Sigma
(X).$ Hence $S_{f_{n_{i}}^{-1}y}=S_{0}$, so $f_{n_{i}}(S_{0})=S_{y},$ for
every $i\in I$. Therefore $y\in f(X).$

Suppose that $\{S_{n},n\in {\Bbb {N} \}}$ has not any constant subnet. By
the compactness of $\Sigma (X)$, there exists a subnet $\{S_{n_{i}},\,i\in
I\}$ of $\{S_{n},n\in {\Bbb {N}\}}$ such that $S_{n_{i}}\rightarrow S$, for
some $S\in \Sigma (X)$. With the above notation, the following is true:

\bigskip

\noindent \underline{{\em Claim.}} {\em There exists a subsequence $%
\{S_{k},\,k\in {\Bbb {N} \}}$ of $\{S_{n},\,\,n\in {\Bbb {N}\}}$ such that
there are $x_{k}\in S_{k}$ with $x_{k}\rightarrow x_{0}$, for some $x_{0}\in
X$. } \newline
{\em Proof.} If not, $R=(\bigcup _{n=1}^{\infty }S_{n})\setminus S$ is
closed in $X$. Indeed, let $R\ni y_{m}\rightarrow y\in X.$ If $y_{m}\in
(\bigcup_{n=1}^{n_{0}}S_{n})\setminus S$ for $m>m_{0}$, then a subsequence
of $\{y_{m},\,m\in {\Bbb {N} \}}$ is contained in some $S_{i}$ for some $%
i\in \{1,\ldots ,n_{0}\}$, therefore $y\in S_{i}\subset R$, as required. If
this is not the case, we construct a subsequence $\{y_{m_{p}},\,p\in {\Bbb {N%
}\}}$ of $\{y_{m},\,m\in {\Bbb {N} \}}$ in the following way: We correspond
to $S_{1}$ the point $y_{m_{1}}\in S_{n_{1}}$ with $n_{1}>1$ and $%
d(y_{m_{1}},y)<1,$ to ($\bigcup _{n=1}^{n_{1}}S_{n})\setminus S$ the point $%
y_{m_{2}}\in S_{n_{2}}$ with $n_{2}>n_{1}$ and $d(y_{m_{2}},y)<\frac{1}{2}$,
and so on. Obviously, $S_{n_{p}}\ni y_{m_{p}}\rightarrow y$, a contradiction.

Hence $R$ is closed in $X$, from which follows $S\subseteq X\setminus R$,
and $X\setminus R$ is open and contains entire components, so $%
S_{n_{i}}\subset X\setminus R,$ eventually. Therefore $S_{n_{i}}=S$, a
contradiction, since we have assumed that $\{S_{n},n\in {\Bbb {N}\}}$ has
not any constant subnet.

\bigskip

According to the Claim, we have $S_{k}\ni x_{k}\rightarrow x_{0}\in X,$
where $S_{k}=S_{f_{k}^{-1}y}=f_{k}^{-1}S_{y},$ from which follows $%
x_{k}=f_{k}^{-1}y_{k}$ for some $y_{k}\in S_{y}.$ Then 
\[
d(y_{k},f(x_{0}))\leq d(y_{k},f_{k}x_{0})+d(f_{k}x_{0},f(x_{0}))=
d(f_{k}^{-1}y_{k},x_{0})+d(f_{k}x_{0},f(x_{0}))\rightarrow 0, 
\]
\noindent therefore $f(x_{0})\in S_{y},$ which means that $S_{y}\cap
f(X)\neq \emptyset $ and, by Lemma 3.5, $S_{y}\subseteq f(X),$ hence $y\in
f(X),$ and $\ f$ is surjective.

\bigskip

\noindent {\large {\bf 3.7}} {\bf Theorem} {\em If $\Sigma (X)$ is compact
then $I(X,d)$ is locally compact.} \newline

\noindent {\bf Proof.} It follows from Lemma 3.3 and Proposition 3.6, since
both conditions 2.1(a) and (b) are satisfied.

\section{The properness of the action $(I(X,d),X)$}

In this short section, applying the methods used previously, we give a
complete proof of the following: \newline

\noindent {\bf Proposition} {\em If $(X,d)$ is locally compact and
connected, then $I(X,d)$ is locally compact and the action $(I(X,d),X)$ is
proper.} \newline

\noindent {\bf Proof.} Since $X$ is connected $G=I(X,d)$ is locally compact
by Theorem 3.7. So, we have to show that, for every $x,$ $y\in X$, there are
neighborhoods $U_{x}$, $U_{y}$ of $x$ and $y$ respectively such that 
\[
(U_{x},U_{y}):=\{g\in G:\,(gU_{x})\cap U_{y}\neq \emptyset \} 
\]
\noindent is relatively compact in $G$. Let $U_{x}=S(x,\varepsilon )$ and $%
U_{y}=S(y,\varepsilon )$ be such that $S(y,2\varepsilon )$ is relatively
compact. Then, for $g\in (U_{x},U_{y})$ and $z\in U_{x}$ with $gz\in U_{y}$,
we have 
\[
d(gx,y)\leq d(gx,gz)+d(gz,y)=d(x,z)+d(gz,y)<2\varepsilon , 
\]
\noindent therefore $g\in F=\{g\in G:\,gx\in S(y,2\varepsilon )\}$. Then $%
x\in K(F)$, and, according to Lemma 3.1, $K(F)$ coincides with the connected
space $X$. From this and Ascoli's theorem it follows that $F$ is relatively
compact in $C(X,X)$. So $(U_{x},U_{y})\subseteq F$ is relatively compact in $%
C(X,X)$, hence in $G$, because $G$ is closed (cf. Proposition 3.6).

\bigskip This proves the Proposition and completes the proof of the Theorem
in the Introduction.

\section{Final Remark}

Using the same arguments we can prove that if $X$ is a locally compact
metrizable space, then $I(X,d)$ is locally compact for all admissible
metrics $d$, provided that the space $Q(X)$ of the quasicomponents of $X$ is
compact (cf. \cite{amanou}). Recall that the quasicomponent of a point is
the intersection of all open and closed sets which contain it. Our
exposition is given via $\Sigma (X)$ because we regard that the condition "$%
\Sigma (X)$ is compact" is a topologically more natural condition than "$Q(X)
$ is compact", althought it is more restrictive: There are locally compact
metric spaces with compact $Q(X)$ and non compact $\Sigma (X)$ as the
following example shows:

\bigskip

\noindent {\bf Example} The space of the connected components of the locally
compact space 
\[
X =\left(\bigcup_{n=1}^{\infty}\{({\frac{1}{n}},y)\;|\,y\in [-1,1]\}\right)
\cup\{(0,y)\;|\,y\in [-1,0)\}\cup \left(\bigcup_{k=1}^{\infty}I_{k}\right)
\subseteq{\Bbb {R}^{2}, }
\]
where 
\[
I_{k}=\{(0,y)\;|\,y\in ({\frac{1}{k+1}},{\frac{1}{k}})\},\;\;\;\; k\in{\Bbb {%
N}^{*}, }
\]
is not compact, because the sequence $\{I_{k}\}\subseteq \Sigma (X)$ does
not have a convergent subsequence in $\Sigma (X)$. On the contrary, $Q(X)$
is compact, because the quasicomponent of the point $(0,-1)$ consists of the
set $\{ (0,y)\;|\; y\in [-1,0)\}$ and the intervals $I_{k}\, ,\; k\in{\Bbb {N%
}^{*}}$.

\bigskip So the compactness of $\Sigma (X)$ is not necessary for the locall
compactness of $I(X,d)$.

\end{document}